\documentclass{amsart} 
\usepackage{amssymb}
\usepackage{amsmath}
\usepackage{amsfonts}

\begin{document}
\newtheorem{theo}{Theorem}[section]
\newtheorem{prop}[theo]{Proposition}
\newtheorem{lemma}[theo]{Lemma}
\newtheorem{exam}[theo]{Example}
\newtheorem{coro}[theo]{Corollary}
\theoremstyle{definition}
\newtheorem{defi}[theo]{Definition}
\newtheorem{rem}[theo]{Remark}


\newcommand{\Bb}{{\bf B}}
\newcommand{\Nb}{{\bf N}}
\newcommand{\Qb}{{\bf Q}}
\newcommand{\Rb}{{\bf R}}
\newcommand{\Zb}{{\bf Z}}
\newcommand{\Ac}{{\mathcal A}}
\newcommand{\Bc}{{\mathcal B}}
\newcommand{\Cc}{{\mathcal C}}
\newcommand{\Dc}{{\mathcal D}}
\newcommand{\Fc}{{\mathcal F}}
\newcommand{\Ic}{{\mathcal I}}
\newcommand{\Jc}{{\mathcal J}}
\newcommand{\Lc}{{\mathcal L}}
\newcommand{\Oc}{{\mathcal O}}
\newcommand{\Pc}{{\mathcal P}}
\newcommand{\Sc}{{\mathcal S}}
\newcommand{\Tc}{{\mathcal T}}
\newcommand{\Uc}{{\mathcal U}}
\newcommand{\Vc}{{\mathcal V}}

\newcommand{\ax}{{\rm ax}}
\newcommand{\Acc}{{\rm Acc}}
\newcommand{\Act}{{\rm Act}}
\newcommand{\ded}{{\rm ded}}
\newcommand{\Gm}{{$\Gamma_0$}}
\newcommand{\ID}{{${\rm ID}_1^i(\Oc)$}}
\newcommand{\PA}{{\rm PA}}
\newcommand{\ACA}{{${\rm ACA}^i$}}
\newcommand{\RefP}{{${\rm Ref}^*({\rm PA}(P))$}}
\newcommand{\RefS}{{${\rm Ref}^*({\rm S}(P))$}}
\newcommand{\Rfn}{{\rm Rfn}}
\newcommand{\tar}{{\rm Tarski}}
\newcommand{\UNFA}{{${\mathcal U}({\rm NFA})$}}

\author{Nik Weaver}

\title [Predicative well-ordering]
       {Predicative well-ordering}

\address {Department of Mathematics\\
          Washington University in Saint Louis\\
          Saint Louis, MO 63130}

\email {nweaver@math.wustl.edu}

\date{\em November 8, 2018}

\begin{abstract}
Confusion over the predicativist conception of well-ordering pervades the
literature and is responsible for widespread fundamental misconceptions about
the nature of predicative reasoning. This short note aims to explain the
core fallacy, first noted in \cite{W}, and some of its consequences.
\end{abstract}

\maketitle


\section{Predicativism}

Predicativism arose in the early 20th century as a response to the
foundational crisis which resulted from the discovery of the classical
paradoxes of naive set theory. It was initially developed in the writings
of Poincar\'{e}, Russell, and Weyl.

Their central concern had to do with the avoidance of definitions they
considered to be circular. Most importantly, they forbade any definition
of a real number which involves quantification over all real numbers.

This version of predicativism is sometimes called ``predicativism given the
natural numbers'' because there is no similar prohibition against defining
a natural number by means of a condition which quantifies over all natural
numbers. That is, one accepts $\mathbb{N}$ as being ``already there'' in
some sense which is sufficient to void any danger of vicious circularity.
In effect, predicativists of this type consider uncountable collections to
be proper classes. On the other hand, they regard countable sets and
constructions as unproblematic.

\section{Second order arithmetic}

Second order arithmetic, in which one has distinct types of variables for
natural numbers ($a$, $b$, $\ldots$) and for sets of natural numbers ($A$,
$B$, $\ldots$), is thus a good setting for predicative reasoning  ---
predicative given the natural numbers, but I will not keep repeating this.
Here it becomes easier to frame the restriction mentioned above in terms of
$\mathcal{P}(\mathbb{N})$, the power set of $\mathbb{N}$, rather than in terms
of $\mathbb{R}$. Thus we reject definitions of sets of numbers of the form
$$\{a \in \mathbb{N}: P(a)\}\eqno{(*)}$$
when $P$ is a formula which includes second order quantifiers, i.e.,
quantification over set variables.

The moral is that predicativists can accept only quite weak comprehension
principles. They cannot reason about $\{a: P(a)\}$ for general $P$.

\section{Countable well-orderings}

Occasionally commentators have worried that predicativist scruples might
forbid any reasoning about uncountable collections. This is wrong: it would
be like saying that a finitist cannot make statements about all natural
numbers, or a platonist cannot make statements about all sets. They can, of
course. In the same way that a finitist can affirm that every natural number
is even or odd, a predicativist can affirm that every real number is rational
or irrational. Indeed, for any particular real number this question can be
answered by a computation of countable length, i.e., it is predicatively
decidable.

An attractive position, expressed first in \cite{W} and later, without credit,
in \cite{F4}, is that in analogy with intuitionists, predicativists should
assume the law of excluded middle only for formulas with no second order
quantifiers. (A formal system embodying this proposal appeared even earlier in
\cite{CP}.) This gives us the ability to reason using such formulas while
neatly explaining the obstruction to forming $\{a: P(a)\}$ as arising from the
possibility that $P(a)$ might not have a definite truth value for some values
of $a$. However, I will not insist on this position here.

In any case, predicativists have a firm grasp of $\mathbb{N}$ and can certainly
affirm that every nonempty set of natural numbers has a least element with
respect to the usual ordering. Given any set of natural numbers, a countable
computation will tell us either that it is empty or what its smallest element
is, and predicativists can see this. So there is nothing impredicative about
defining a total ordering $\preceq$ of $\mathbb{N}$ to be a well-ordering
if every nonempty subset of $\mathbb{N}$ has a $\preceq$-least element.
Yes, this definition involves quantification over $\mathcal{P}(\mathbb{N})$,
but there is no circularity issue because it does not introduce any new
objects. It does nothing more than fix the term ``well-ordered'' as an
abbreviation of the longer, predicatively intelligible expression ``every
nonempty subset has a least element''.

\section{Induction for sets}

The following formulation of the well-ordering property is also useful. Say
that a subset $A \subseteq \mathbb{N}$ is {\it progressive} for a total
order $\preceq$ on $\mathbb{N}$ if $(\forall b \prec a)(b \in A)$ implies
$a \in A$. That is, a number belongs to $A$ whenever every preceding number
belongs to $A$. Now consider the condition
\begin{description}
\item[(A)] If $A \subseteq \mathbb{N}$ is progressive for $\preceq$ then
$A = \mathbb{N}$.
\end{description}
This is equivalent to saying that $\preceq$ is a well-order: in one direction,
if $\preceq$ is a well-order and $A$ is progressive, then $\mathbb{N}
\setminus A$ cannot have a least element by progressivity, and hence it must
be empty; conversely, if condition (A) holds and $A \subseteq \mathbb{N}$ is
nonempty, then $\mathbb{N}\setminus A$ cannot be progressive, so there must
exist $a \in A$ such that $(\forall b \prec a)(b \not\in A)$, i.e., $a$ is the
least element of $A$.

I have detailed this simple argument because I want to make the point that
this equivalence is predicatively valid. That is, the argument which shows
that condition (A) is equivalent to $\preceq$ being a well-order is one which
a predicativist could make. It should be easy to convince oneself of this.

\section{Induction for predicates}

Say that a predicate $P$ --- a formula with one free number variable --- is
progressive for $\preceq$ if $(\forall b \prec a)P(b)$ implies $P(a)$. We
have the following stronger version of condition (A).
\begin{description}
\item[(B)] If a predicate $P$ is progressive for $\preceq$ then
$(\forall a \in \mathbb{N})P(a)$.
\end{description}

This condition is stronger in the sense that it reduces to condition (A)
when we take $P(a)$ to be the predicate ``$a \in A$''. However, condition (B)
can also be inferred from condition (A) by the following simple argument:
\begin{enumerate}
\item Assume condition (A) and suppose $P$ is progressive for $\preceq$.

\item Let $A = \{a \in \mathbb{N}: P(a)\}$.

\item Since $P$ is a progressive predicate, $A$ is a progressive subset.
Thus condition (A) yields $A = \mathbb{N}$, which is just to say that $P(a)$
holds for all $a \in \mathbb{N}$.
\end{enumerate}
But this is an argument predicativists cannot generally make. Step 2 invokes a
comprehension principle which is impredicative unless $P$ has a special form.

Condition (B) does not predicatively follow from condition (A). This is what no
one has understood.

\section{But doesn't it really?}

Most working mathematicians have probably internalized the equivalence of
conditions (A) and (B), making it seem so obvious as not to need any special
argument. But the two are not a priori equivalent. Condition (B) is stronger.
It is only by invoking a comprehension axiom that one gets from (A) to (B).

Some people I have spoken to have acknowledged that (A) and (B) are a priori
inequivalent, but have maintained that a predicativist can somehow just
intuit that (B) follows from (A). Again, I think this point of view arises from
having internalized the equivalence. The two statements are not the same. To
get from one to the other you need to invoke a comprehension principle, and
predicativists cannot do this. If predicativism means anything, it entails
weak comprehension axioms. We cannot simply postulate that predicativists
can ``just see'' the validity of some argument that relies on impredicative
principles. One could draw any conclusion whatever from that kind of reasoning.

A more subtle idea is that it may in fact be the case that whenever there is
a predicatively valid argument which shows that $\preceq$ is well-ordered in
the sense of condition (A), there will be some other predicatively valid
argument which shows that it verifies condition (B), at least schematically.
But I do not see how one could prove this without first delineating precisely
which ordinals are predicatively available.

\section{Feferman-Sch\"{u}tte}

The celebrated Feferman-Sch\"{u}tte analysis \cite{F1, S} of predicatively
provable ordinals allegedly provides just such a delineation. But it is wrong.

The goal is to determine which ordinals are isomorphic to recursive total
orderings of $\mathbb{N}$ which can be proved to be well-orderings by
predicative means. The Feferman-Sch\"{u}tte analysis can be presented in
several ways, but the most intuitive involves a system which allows proof
trees of infinite (well-ordered) heights over a straightforwardly predicative
base system. The idea is that one starts with proof trees of
height $\epsilon_0$, say, but once one has proven that a notation for
$\alpha$ is well-ordered --- so that $\alpha$ is predicatively ``recognized''
to be an ordinal --- one is allowed to use proof trees of height $\alpha$.

If we define $\gamma_1 = \epsilon_0$ and $\gamma_{n+1} = \phi_{\gamma_n}(0)$,
where $\phi_\alpha$ are the Veblen functions, then using a proof tree of
height $\gamma_n$ one can prove that a notation for $\gamma_{n+1}$ is
well-ordered. Thus, iterating this process, every ordinal less than $\Gamma_0 =
\sup \gamma_n$ is
supposed to be predicatively provable. Conversely, one cannot prove that a
notation for $\Gamma_0$ is well-ordered using a proof tree of any height
less than $\Gamma_0$, and this allegedly shows that $\Gamma_0$ is the limit
of predicative reasoning.

An obvious question is why the predicativist cannot recognize for himself
that for all $n$ he can prove a notation for $\gamma_n$ is well-ordered, and
infer from this that $\Gamma_0$ is well-ordered. It seems like just the sort
of countable reasoning that predicativists are good at. This question was
raised in \cite{H}, among other places. I will return to it in a moment.

But first, let us ask how the predicativist is to infer the soundness of proof
trees of height $\alpha$ from the fact that $\alpha$ is well-ordered. Why
should it follow that such trees prove true theorems? It is easy to see that
soundness is progressive --- if every proof tree of height less than $\alpha$
is sound, then every proof tree of height $\alpha$ is sound --- so one wants
to use an induction argument in the form of condition (B) with $P(a) =$ ``proof
trees of height $\alpha$ are sound'', where $a$ is a notation for $\alpha$.

But all we prove with a tree of height $\gamma_n$ is that a notation for
$\gamma_{n+1}$ is well-ordered in the sense of condition (A). Lacking an
impredicative comprehension principle, we cannot infer the soundness of
proof trees of height $\gamma_{n+1}$.

The Feferman-Sch\"{u}tte analysis relies essentially on an impredicative
inference of condition (B) from condition (A).

\section{Induction vs.\ recursion}

Now Feferman-Sch\"{u}tte does not really need the soundness of all proof
trees of height $\alpha$. Closer examination of the argument reveals that
what is essential is to be able to carry out recursive constructions of
length $\alpha$. That is, given that for any $\beta \prec \alpha$
and any indexed family $(A_\gamma)_{\gamma \prec \beta}$ there exists a
unique $A_{\beta} \subseteq \mathbb{N}$ which satisfies some condition
relative to the previous $A_{\gamma}$'s, we need to conclude that there
is a family $(A_\gamma)_{\gamma \prec \alpha}$ which satisfies the condition
at each stage. That is all we really need. But this still requires an instance
of condition (B) with $P$ containing one second order quantifier.

Thus, even stripped to its bare minimum, the argument which gets one from
$\gamma_n$ to $\gamma_{n+1}$ requires impredicative comprehension.

A prominent logician recently told me that ``everyone agrees'' with the
Feferman-Sch\"{u}tte analysis. I think this is basically true. Because
no one has appreciated the distinction between induction for sets and
induction for predicates.

One possible exception is Feferman himself: buried in his little-known paper
\cite{F2} is a comment that ``the well-ordering statement $\ldots$ on the
face of it only impredicatively justifies the transfinite iteration of
accepted principles up to $a$,'' which looks a lot like a direct
acknowledgement of the defect we have been discussing. But then in
the later paper \cite{F3} he refers to ``prima facie impredicative notions
such as those of ordinals or well-orderings,'' revealing a basic
misunderstanding of predicativism (cf.\ Section 3 above), and in \cite{F3}
he also affirms that a version of the inference from condition (A) to
condition (B) is
predicative because it ``accords with ordinary informal reasoning''.

In a personal communication to me, Feferman firmly denied that he had {\it not}
repudiated the autonomous systems analysis, which I suppose means that he
felt he had in fact repudiated it. His position, I think, was that later papers
of his such as \cite{F2} or \cite{F3} remedied this defect by analyzing the
problem in a different way. However, these later systems are shown in
\cite{W} to suffer from similar defects (and indeed, were obviously
deliberately engineered to get the answer $\Gamma_0$).

\section{A wild hope}

It may be hard for anyone who has already bought into the Feferman-Sch\"{u}tte
analysis to accept that it is wrong, no matter how plain the error is. But I
do not expect anyone to attempt the futile task of trying to find some
predicatively valid argument that infers recursion up to $\alpha$ from
induction up to $\alpha$. The lazy response, the only one I have personally
encountered, is to simply postulate that a predicativist can directly intuit
this inference without needing to prove it.

Besides being absurd on its face, this response runs into the obvious question
mentioned in Section 7, about the predicativist being able to recognize
that for all $n$ he can prove a notation for $\gamma_n$ is well-ordered. He
would then be able to get beyond $\Gamma_0$. So if one really wants to stop
exactly at $\Gamma_0$, one has to postulate that the predicativist is able to
magically intuit the inference of recursion from induction in any particular
case, {\it but he cannot recognize that he has this general ability}, as this
would get him past $\Gamma_0$. For every $n$, when he proves that a notation
for $\gamma_n$ is well-ordered, the revelatory intuition that it also supports
recursion has to come as a surprise. That is the depth of irrationality to
which one must sink in order to preserve $\Gamma_0$ as the limit of
predicative reasoning.

Let me return to the more subtle idea mentioned in Section 6. Could it be
that whenever a predicativist can prove that a total order is a
well-order, he can also prove in some different way that it supports
transfinite recursion? As I said before, I do not see how one could prove
this in general, but proving it for all ordinals less than $\Gamma_0$, say,
does not seem out of the question.

In fact, I showed in \cite{W} that using hierarchies of Tarskian truth
predicates one can prove stronger forms of condition (B) at lower ordinals
which allow one to prove weaker forms of (B) at higher ordinals, sufficiently
to get up to $\Gamma_0$. In my opinion these techniques are predicatively
legitimate. The problem for Feferman-Sch\"{u}tte is that using them one can
continue on well past $\Gamma_0$, as was shown in detail in \cite{W}.

\section{Inductive definitions}

The Feferman-Sch\"{u}tte analysis is not the only place where the fallacy of
conflating conditions (A) and (B) appears. For instance, consider the theory
$ID_1$ of one inductive definition. It has been called ``generalized
predicative'', although according to the Stanford Encyclopedia of Philosophy
it is simply ``considered predicative in today's foundations of constructive
mathematics'' \cite{C}. Its proof theoretic ordinal is much larger than
$\Gamma_0$.

In $ID_1$ one wants to introduce a class as the smallest class contained in
$\mathbb{N}$ which respects some closure operation. One formalizes this
idea by introducing a constant symbol $\mathcal{C}$ to represent the class
and expressing its minimality by an axiom scheme which states, for any
predicate $P$, that if $P$ respects the closure condition then every
element of $\mathcal{C}$ satisfies $P$.

Classically $\mathcal{C}$ exists as the intersection of all classes which
respect the closure condition, but this is not a predicatively valid
definition: it defines a set of numbers by a condition which quantifies over
all sets of numbers. However, it seems like one ought to be able to
predicatively regard $\mathcal{C}$ in roughly the same way one predicatively
regards $\mathbb{R}$, as a proper class which can be built up in stages but
never reaches a finished state.

This is not right. The suggestion that $ID_1$ is predicative in any sense is
wrong, because the minimality
axioms are impredicative. To see this, imagine how one would establish that
$P(a)$ holds for every $a \in \mathcal{C}$, given that $P$ respects the closure
operation. This hypothesis tells us that if every element appearing at some
stage of the construction of $\mathcal{C}$ satisfies $P$, then every element
at the next stage will too. And at limit stages one merely collects together
all the preceding stages, so nothing new appears and thus the minimality axiom
is still verified, provided it was satisfied at all previous stages.

In other words, the condition ``every element in the construction of
$\mathcal{C}$ at stage $\alpha$ satisfies $P$'' is progressive in $\alpha$.
So if one imagines constructing $\mathcal{C}$ in a well-ordered series of
stages, then the minimality axiom for $P$ should be satisfied --- provided
we can infer condition (B) from condition (A). Which we cannot predicatively
do. It is the same error that invalidates Feferman-Sch\"{u}tte.

Another idea is to specify that $\mathcal{C}$ is to be constructed in stages
which are well-ordered not in the sense of condition (A), but in the stronger
sense of supporting induction for any predicate of the language. But no, this
is badly impredicative because the language contains the constant symbol
``$\mathcal{C}$'', which has no meaning until one has specified what
$\mathcal{C}$ is, so that one would be explaining how $\mathcal{C}$ is to
be constructed in terms of conditions which refer to $\mathcal{C}$.

\section{Kripke-Platek and CZF}

Kripke-Platek set theory (KP) and Constructive Zermelo-Fraenkel set theory
(CZF) are two set theoretic systems which are also routinely claimed to be
predicative. (According to Wikipedia, KP is ``roughly the predicative part
of ZFC'' and CZF has ``a fairly direct constructive and predicative
justification''.)

In fact, both are impredicative for the same reason $ID_1$ is: yet again, the
fallacy involves a confusion between conditions (A) and (B). In both cases
the problematic axioms are the set induction scheme, which states, for any
formula $P$,
$$(\forall y)([\forall x \in y\, P(x)] \to P(y)) \to (\forall y)P(y).$$
Informally, if a predicate holds of a set $y$ whenever it holds of all the
elements of $y$, that predicate must hold of all sets.

The informal justification for this scheme hinges on the premise that the
universe of sets is built up in a well-ordered series of stages. One then
applies progressivity of $P$ to infer, inductively, that it holds of all
sets in the universe.

Just as with $ID_1$, this justification fails because being well-ordered
in the sense of condition (A) does not predicatively entail the instances
of condition (B) which would be needed to make the induction argument.
And also as in that case, there is no option of strengthening the premise to
say that the universe of sets is built up in a series of stages which are
well-ordered in some stronger way which affirms condition (B). The instances
of condition (B) which we would need in order to justify set induction
involve all predicates expressible in the language of set theory, but the
latter does not have an interpretation until we specify how the universe of
sets is to be built up. So this would be circular.

\section{Future directions}

All this was explained in great detail in \cite{W}. As I mentioned above,
an attempt was also made there to rehabilitate the argument that gets us
up to $\Gamma_0$ by means of iterated truth predicates. But when one does
this there is no particular reason to stop at $\Gamma_0$; the iterated truth
technique gets us up to at least the small Veblen ordinal
$\phi_{\Omega^\omega}(0)$. However, the farther one goes the more complicated
the analysis becomes.

Very likely the argument given in \cite{W} can be simplified. I would still
expect it to become more complicated the higher one goes. In my opinion
predicativists can get past $\Gamma_0$, but not up to the Bachmann-Howard
ordinal --- the evidence for this is that, as discussed above, all known
systems which get that far are clearly impredicative. The exact ordinal
limit of predicative reasoning might not be a well-defined concept.

\end{document}